\documentclass[english]{article}
\usepackage{amscd,amssymb,stmaryrd}
\usepackage[fleqn]{amsmath}
\usepackage{graphicx}
\usepackage{caption}

\usepackage[utf8]{inputenc}
\usepackage[export]{adjustbox}
\usepackage{wrapfig}
\usepackage{amstext}
\usepackage{pstricks,pst-node,pst-plot,pst-coil}
\usepackage{amsthm}
\usepackage{amsmath}
\usepackage{color}
\usepackage{mathtools}
\usepackage{hyperref}
\usepackage[english]{babel}
\usepackage{booktabs}
\usepackage{mathrsfs}
\usepackage{multirow}
\usepackage{subfiles}
\usepackage{subfigure}
\usepackage{lineno}
\usepackage{floatrow}
\usepackage{authblk}

\newfloatcommand{capbtabbox}{table}[][\FBwidth]

\providecommand{\keywords}[1]
{
  \small	
  \textbf{\textit{Keywords---}} #1
}
\title{Online conservative generalized multiscale finite element method for flow models}
 \author[1]{Yiran Wang}
 \author[1]{Eric Chung}
\author[2]{Shubin Fu}
\author[3]{Michael Presho}
\affil[1]{Department of Mathematics, The Chinese University of Hong Kong, Hong Kong SAR}
\affil[2]{Department of Mathematics, University of Wisconsin-Madison,WI, USA}
\affil[3]{Department of Mathematics, Southeast Missouri State University, Cape Girardeau, MO 63701, United States}
\begin{document}

\maketitle

\begin{abstract}
In this paper, we consider an online enrichment procedure using the Generalized Multiscale Finite Element Method (GMsFEM) in the context of a two-phase flow model in heterogeneous porous media. The coefficient of the elliptic equation is referred to as the permeability and is the main source of heterogeneity within the model. The elliptic pressure equation is solved using online GMsFEM, and is coupled with a hyperbolic transport equation where local conservation of mass is necessary. To satisfy the conservation property, we aim at constructing conservative fluxes within the space of multiscale basis functions through the use of a postprocessing technique. In order to improve the accuracy of the pressure and velocity solutions in the online GMsFEM we apply a systematic online enrichment procedure. The increase in pressure accuracy due to the online construction is inherited by the conservative flux fields and the desired saturation solutions from the coupled transport equation. Despite the fact that the coefficient of the pressure equation is dependent on the saturation which may vary in time, we may construct an approximation space using the initial coefficient where no further basis updates follow. Numerical results corresponding to four different types of heterogeneous permeability coefficients are exhibited to test the proposed methodology.
\end{abstract}


\section{Introduction}
A large number of problems of fundamental and practical significance are described by partial differential equations with coefficients that vary over a wide range of length scales. For example, composite materials, porous media and turbulent transport in high Reynolds number flows are models of this type. The heterogeneity and high-contrast properties of the coefficients cause significant difficulty in analyzing these types problems. In this paper, we consider a two-phase flow model in which the so-called permeability coefficient is assumed to be highly heterogeneous. Solving this type of model problem on a fine scale that sufficiently captures the underlying behavior of the heterogeneity may become prohibitively expensive. As a result, methods that aim toward effectively reducing the dimension of the associated fine-scale system(s) have been a topic of continued interest in recent decades. For example, upscaling procedures (see, e.g., \cite{chen2003coupled,durlofsky1991numerical,wu2002analysis}) and multiscale methods (see, e.g., \cite{efendiev2011multiscale,jenny2003multi,wheeler2012multiscale,efendiev2009multiscale,hou1997multiscale} ) are approaches that have been shown to offer effective alternatives to direct fine-scale computations.
For upscaling, one derives a set of localized problems in which averaged quantities may be maintained while solving a lower dimensional global problem on a coarse grid. However, this type of approach may diminish important fine-scale information that has strong effect on the solution behavior. Multiscale methods, on the other hand, hinge on the the independent construction of a set of multiscale basis functions that are used to span a coarse-grid solution space. The coarse-grid discretization parameter may be much larger than the characteristic scale of heterogeneous coefficient, however, the multiscale basis functions inherently include the fine-scale information of the underlying heterogeneity of the medium.

In order to model multi-phase flow, local mass conservation for the fluid velocity fields is required. This requirement has motivated a variety of mass-conservative approaches, such as multiscale finite volume methods \cite{cortinovis2014iterative,jenny2003multi,lunati2004multi}, mixed multiscale finite element methods \cite{aarnes2004use,aarnes2008mixed,chen2003mixed,chung2015mixed}, mortar multiscale methods \cite{arbogast2007multiscale, peszynska2005mortar,peszynska2002mortar}, discontinous Galerkin (DG) methods \cite{du2018adaptive,kim2013staggered,cockburn2002local}, and postprocessing methods \cite{odsaeter2017postprocessing,bush2013application}. In this paper, we use a global continuous Galerkin (CG) method. An advantage of CG multiscale formulation is the relative ease of implementation. On the other hand, a CG solution does not automatically satisfy local conservation, which is essential in our model problem. In order to address this limitation, we adpot an analogous postprocessing technique from \cite{bush2014application}. In particular, after obtaining a multiscale solution, we solve an independent set of local auxiliary problems in order to obtain the locally conservative fluxes.

In terms of the standard Multiscale Finite Element Method (MsFEM), there are two main shortcomings. The first one is that only one basis in each local neighborhood may not be sufficient to guarantee an accurate approximation, especially when there are long channels and non-separable scales in the permeability field. The second limitation is that we often assume that the local boundary conditions are linear along the edges of coarse blocks, which may create a mismatch between multiscale solution and fine-scale solution on the coarse block boundaries. One such technique that may be used to reduce the effect of boundary terms is oversampling \cite{efendiev2000convergence,hou1997multiscale}. Oversampling involves the enlargement of the local computing regions in order to address the linear boundary values. A more recent method that serves to improve the accuracy of MsFEM is the Generalized Multiscale Finite Element Method (GMsFEM) \cite{efendiev2013generalized}.
GMsFEM is a flexible general framework that generalizes MsFEM by systematically enriching the coarse spaces. In particular, more basis functions are added to the initial approximation space in order to improve the accuracy of the multiscale solution. The creation of GMsFEM solution spaces often involves the construction of snapshot, offline, and online spaces \cite{chung2015residual,chung2017online} in order to streamline the procedure for repeated basis function computations. In order to construct an offline multiscale space, some well-designed local spectral problems are solved in order to obtain a set of basis functions that are independent of global information such as source terms and boundary conditions. These local problems are motivated by the convergence analysis, which offers a convergence rate of $1/\Lambda$, where $\Lambda$ is the smallest eigenvalue whose modes are excluded in the multiscale space. If we increase the number of offline bases to a certain number, the error decay will diminish, and in \cite{CHUNG201669}, it is shown that a good approximation from the reduced model can be expected only if the offline information is a good representation of the problem.  Consequently, an online enrichment procedure is essential if the offline bases are not sufficiently accurate. The main idea in this paper is to enrich the offline space by incorporating a new set of basis functions in order to obtain a significant error decay. In \cite{chan2016adaptive,chung2015residual}, the authors propose an online construction resulting from the associated offline space. In consideration of fact that the offline bases are obtained independently through a set of local problems, one may seek to construct a set of bases that contain some global information. Based on this idea, we use residual-driven basis functions which are computed through a set of local problems. The analysis in \cite{chung2017online,chung2015residual} shows that the error decay is proportional to $1-\Lambda$.

The rest of paper is organized as follows. In Section 2 we introduce the model problem and the corresponding solution algorithm. In Section 3, we describe the Generalized Multiscale Finte Element Method (GMsFEM) and the construction of the online solution space. The post-processing technique that is used to ensure the local conservation of mass property is reviewed in Section 4. In Section 5 we offer a variety of numerical results to illustrate the effectiveness of the proposed methodology.

\section{Model problem}

\subsection{Two-phase model}
In this paper, we consider the dynamics of  the movements of two immiscible fluids in a heterogeneous oil reservoir constrained in a domain $\Omega$. In particular, we model scenario where water is discharged to replace trapped oil in a saturated subsurface. Under the assumptions that the environment is gravity-free, capillary pressure is not included, and that two fluids fill the pore space we can apply Darcy's law combined with a statement of conservation of mass. The principle equations of the flow may then be stated as follows:
\begin{eqnarray}
\nabla\cdot \boldsymbol{v}=q, ~~ \text{where} ~~ \boldsymbol{v}=-\lambda(S)k(x)\nabla p\label{elliptic}
\end{eqnarray}
\begin{eqnarray}
\frac{\partial S}{\partial t}+\nabla \cdot(f(S) \boldsymbol{v})=q_{w}, \label{transport}
\end{eqnarray}
where $p$ is the pressure, $\mathbf{v}$ is the Darcy velocity, $S$ is the water saturation, $q,q_w$ are any external forces and $k(x)$ is the heterogeneous permeability coefficient. The total mobility $\lambda(S)$ and the flux function $f(S)$ are respectively given by:
\begin{eqnarray*}
\lambda(S)=\frac{k_{r w}(S)}{\mu_{w}}+\frac{k_{r o}(S)}{u_{e}}, \quad f(S)=\frac{k_{r w}(S) / \mu_{w}}{\lambda(S)}
\end{eqnarray*}
where $k_{r,j},j=w,o$, is the relative permeability of the phase $j$.
\subsection{Solution algorithm}
In Table \ref{two-phase algorithm}, we display the algorithm that is used to solve the two-phase model in Eqs. \eqref{elliptic} and \eqref{transport}. In order to solve for the unknown saturation $S$, we first split the time interval into a set of specified subintervals. $S$ is initialized by $S_0$ and then solved by a series of iterations that are included in Table \ref{two-phase algorithm}. More specifically, we use $S_{n-1}$ in \eqref{elliptic} to obtain $p_n$ and $v_n$. Then we solve \eqref{transport} using the new flux $\mathbf{v}_n$ to obtain $S_n$.\\
\begin{table}[!htbp]
    \begin{tabular}{c l}
    \hline
      &Two-phase algorithm  \\
       \hline
    \textbf{Input} & $S_{n-1}$ obtained in previous time step\\
      \textbf{Output}& $S_{n}$ \\
      &1. Solving \ref{elliptic} to get $p_n$ and $\mathbf{v}_n$\\
      &2. Using $\boldsymbol{v}_n$ and $S_{n-1}$  in \ref{transport} to get $S_{n}$ \\
      \hline
      \caption{Two-phase algorithm}
      \label{two-phase algorithm}
      \end{tabular}
\end{table}

To solve \eqref{transport}, we integrate over the time interval $[t_{n-1},t_n]$ and a control volume  $C_z\subset \Omega$ to obtain
\begin{eqnarray}
\text{meas}(C_z)(S_{z,n}-S_{z,n-1})+\Delta t\int_{\partial_{C_z}}\mathbf{v}\cdot \mathbf{n} f(S_{z,n-1}) ~ dl=\Delta t\int_{C_z} q_w ~ dx,
\end{eqnarray}
where we have neglected the error terms, and we use
\begin{eqnarray}
S_{z,n}\approx\dfrac{1}{meas(C_z)}\int_{C_z}S(x,t_n) ~ dx.
\end{eqnarray}
We use $meas(A)=\int_{\Omega} 1_{A}d x$ with $1_A=1$ when $x\in A$ while 0 elsewhere.
To evaluate the term $\int_{\partial_{C_z}} \mathbf{v} \cdot \mathbf{n} f(S_{z,n-1}) \, dl$, we use an upwinding scheme. A review of upwinding on a rectangular mesh can be in \cite{thomas2013numerical}, for example. It is imperative that the numerical approximation of $\mathbf{v}$ satisfies the following local conservation property. In particular, it is desirable to have
\begin{eqnarray}
\int_{\partial{C_z}}\mathbf{v}\cdot \mathbf{n} ~d l=\int_{C_z} q ~d x.
\end{eqnarray}
There are two main ways to obtain the desired quantities $(\mathbf{v},p)$. The first one is to simultaneously solve the first order system \eqref{elliptic}. For example, one may apply the mixed finite element formulation \cite{chan2016adaptive}. In this paper, we consider the alternative of transforming \eqref{elliptic} into a second order equation that  governs the pressure  $p$. The approximation of $\mathbf{v}$ is calculated using the relation $\mathbf{v}=-\lambda(S)k(x)\nabla p$, and a postprocessing procedure follows for local conservation. Since it is computationally expensive to apply the postprocessing procedure on the fine-scale solution, we instead use the Generalized Multiscale Finite Element Method (GMsFEM), which will be introduced in the next section.

\section{Generalized multiscale finite element method}
\subsection{Preliminaries}
\begin{figure}[!htbp]
\centering
\subfigure[$\kappa_1(x)$]
{ \includegraphics[width=0.45\textwidth]{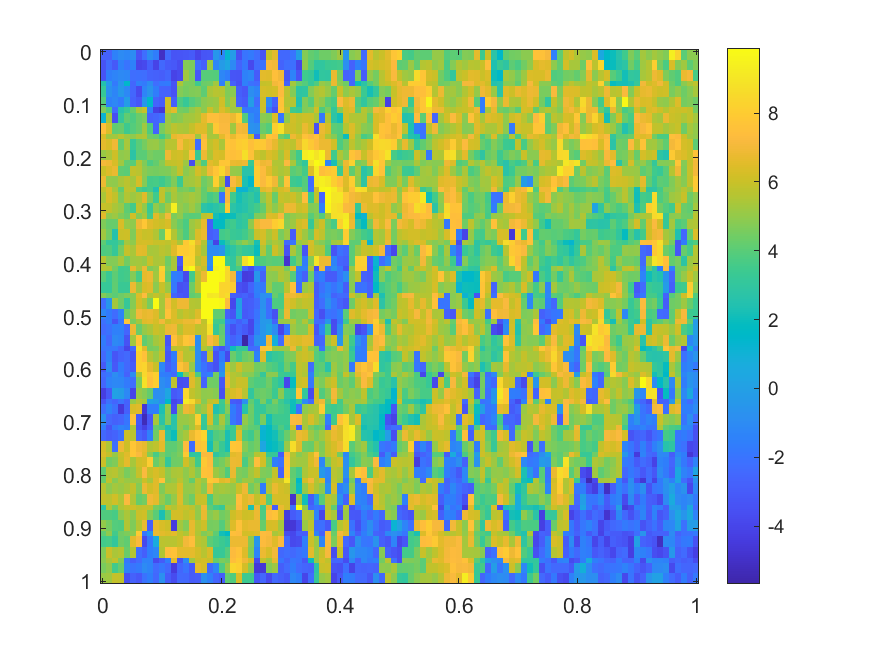}}
\subfigure[$\kappa_2(x)$]
{ \includegraphics[width=0.45\textwidth]{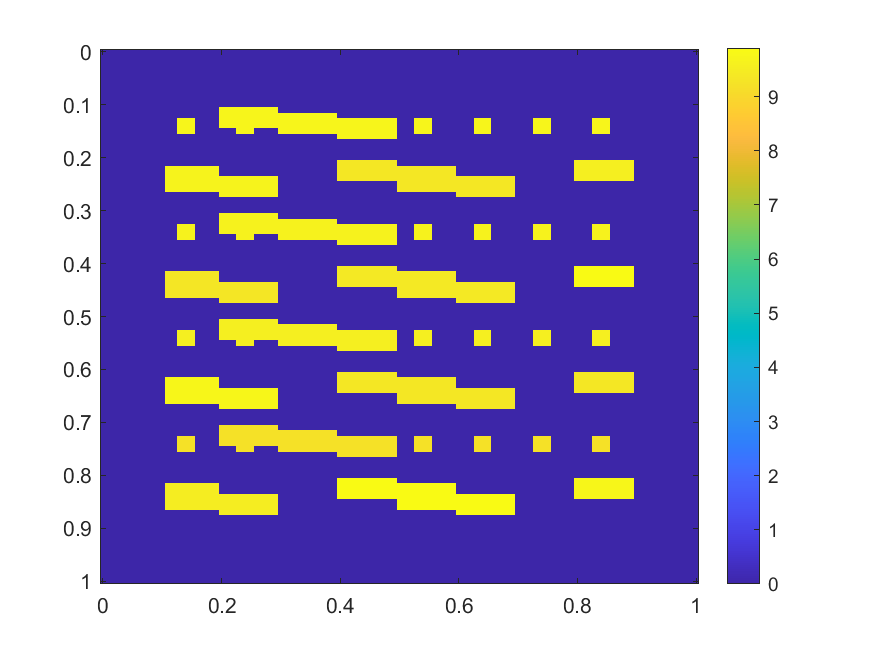}}
\subfigure[$\kappa_3(x)$]
{ \includegraphics[width=0.45\textwidth]{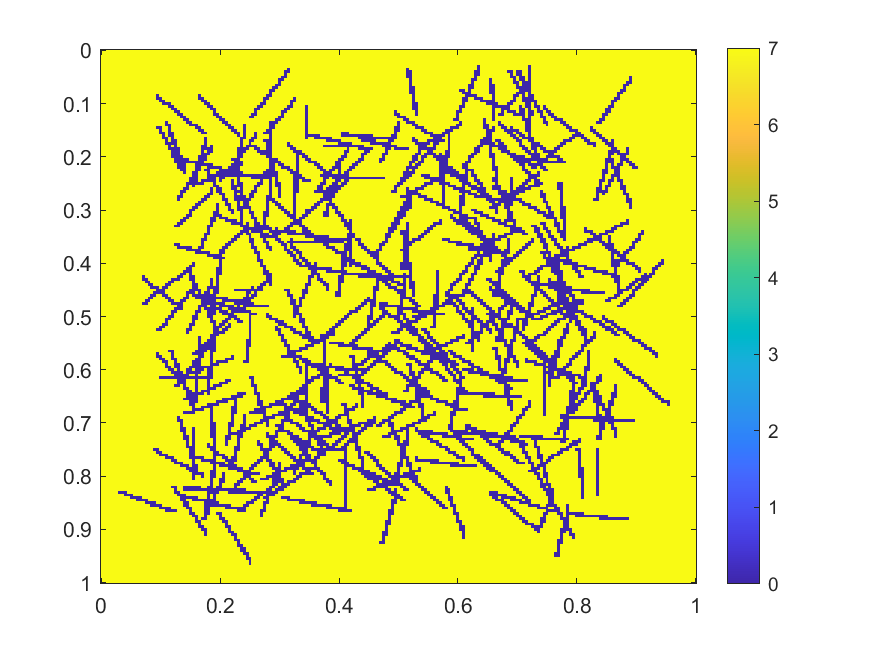}}
\subfigure[$\kappa_4(x)$]
{ \includegraphics[width=0.45\textwidth]{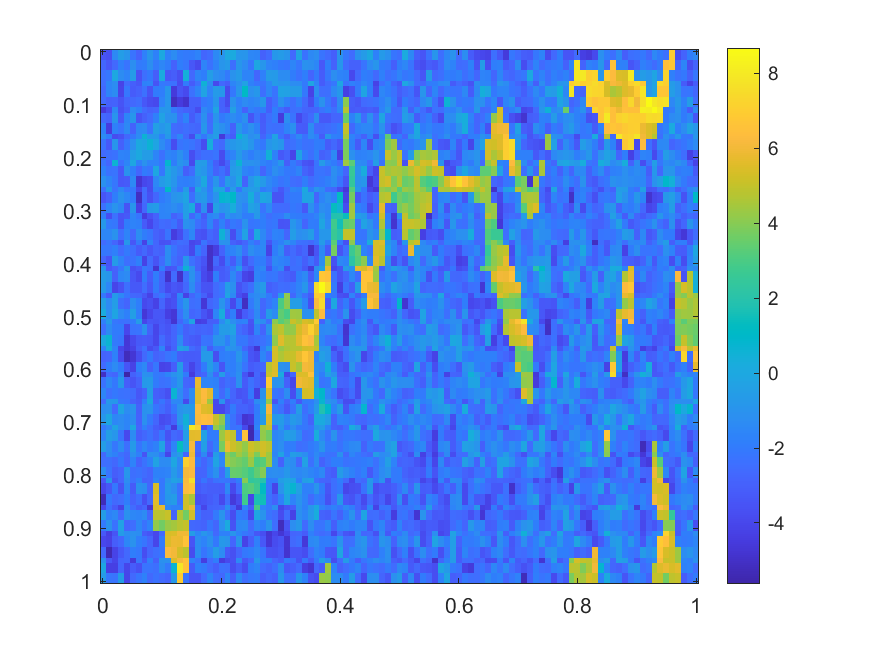}}
\caption{Examples of heterogeneous permeability fields; all plots are on the log scale}\label{medium}
\end{figure}

We fix our attention to the following second order elliptic problem
\begin{eqnarray}
\begin{aligned}
-\operatorname{div}(\lambda k(x) \nabla p)&=q \quad \text { in } \Omega \\
p&=p_{D} \quad \text { on } \Gamma_{D} \\
-\lambda k \nabla p \cdot \mathbf{n}&=g_{N} \quad \text { on } \Gamma_{N}\label{model}
\end{aligned}
\end{eqnarray}
where $k(x)$ is a highly heterogeneous field with high contrast. In practice, we assume that there is a positive  constant $k_{min}$ such that $k(x)\geq k_{min}\geq0$, while $k(x)$ can vary widely (i.e.,
$k_{max}/k_{min}$ is very large, for example $10^5$). Four examples of $k(x)$ that are considered in this paper are offered in Figure \ref{medium}. All permeability fields in the figure are plotted on the log scale. Additionally, $\lambda$ is a known mobility coefficient, $q$ denotes any external forcing, and $p$ is an unknown pressure field satisfying Dirichlet and Neumann boundary conditions given by $p_D$ and $g_N$, respectively. Here $\Omega$ is a convex polygonal and two dimensional domain with boundary $\partial{\Omega}=\Gamma_{D}\cup \Gamma_{N}$.

We consider a function in $H^1(\Omega)$ whose trace on $\Gamma_{D}$ coincides with the given value $p_D$; we
denote this function also by $p_D$. The variational formulation of \eqref{model} is stated as follows. We find $p\in H^1({\Omega})$ with $(p-p_D)\in H^1_{\Omega}=\{w\in H^1_0(\Omega):w|_{\Gamma_{D}}=0 \}$ such that
\begin{eqnarray}
a(p,v)=F(v)-\left\langle g_N,v\right\rangle_{\Gamma_{N}} ~~\text{for all } ~~v\in H^1_D \label{varational}
\end{eqnarray}
where
\begin{eqnarray*}
\begin{aligned}
a(p, v)&=\int_{\Omega} \lambda k(x) \nabla p(x) \nabla v(x)~ d x,\\
F(v)&=\int_{\Omega}q(x)v(x) ~d x, ~~\text{and}\\
\left\langle g_{N}, v\right\rangle_{\Gamma_{N}}&=\int_{\Gamma_{N}} g_{N}(x) v(x)~ d l.
\end{aligned}
\end{eqnarray*}
 \begin{figure}[ht]
    \centering
 {\includegraphics[width=3in]{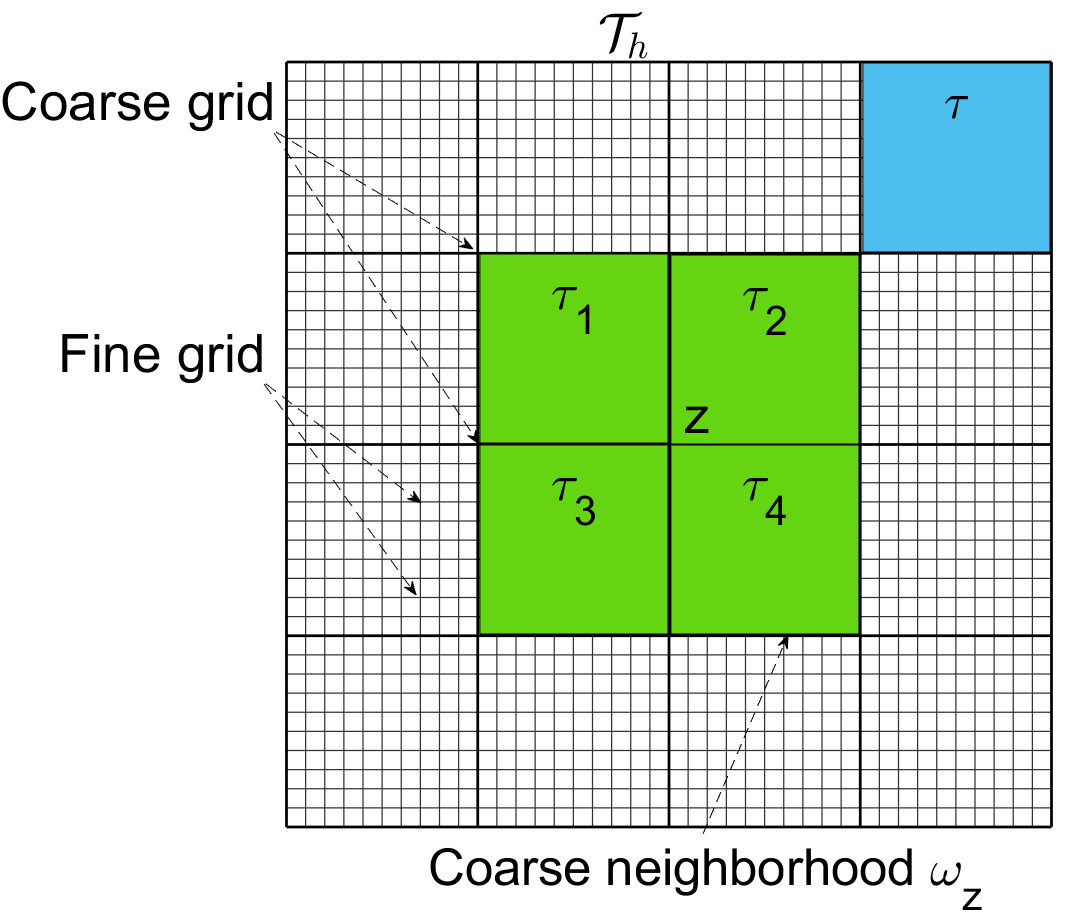}}
 \caption{Discretization of $\Omega$ into $\mathcal{T}_h=\cup \tau$. Here $\omega_z=\cup_{i=1}^4\tau_i$ is the supp($\chi_z$).}\label{figure:coarse}
 \end{figure}

In order to implement a finite element approximation of \eqref{varational}, we
let $\mathcal{T}^{h}$ denote a partition of the domain $\Omega$ into fine elements. Here, $h>0$ is used to denote the fine-grid mesh size.
 The coarse partition, $\mathcal{T}^{H}$ of the domain $\Omega$, is formed such that each element in $\mathcal{T}^{H}$ is a connected union of fine-grid blocks. More precisely, $\forall K_{j} \in \mathcal{T}^{H}$, $ K_{j}=\bigcup_{F\in I_{j} }F$ for some $I_{j}\subset \mathcal{T}^{h}$. The quantity $H>0$ is the coarse mesh size. In this paper we consider the case of rectangular coarse elements, yet the methodology can be used with general coarse elements. An illustration of the mesh notations is shown in the Figure \ref{figure:coarse} (\textbf{the notation in the illustration does not match the notation used below. For example $\omega_z$ is used in the figure for a neighborhood, whereas $D_i$ is used below for the neighborhood}). We denote the interior nodes of $\mathcal{T}^{H}$ by $x_i, ~~i=1,\cdots,N_{\text{in}}$,
 where $N_\text{in}$ is the number of interior nodes. The coarse elements
 of $\mathcal{T}^{H}$ are denoted by $K_j, ~~j=1,2,\cdots,N_e$, where $N_e$ is the number of coarse elements. We define the coarse neighborhood of the nodes $x_i$ by $D_i:=\cup\{K_j\in T^{H}:x_i\in \overline{K_j}\}$.

\subsection{GMsFEM for pressure equation}
 In this paper, we will apply the GMsFEM to solve nonlinear parabolic equations. The
 method is motivated by the finite element framework. First, a variational formulation is defined. Then we construct some multiscale basis functions.
  Once the fine grids are given, we can compute the fine-grid solution. Let $\gamma_1,\cdots,\gamma_n$ be the standard finite element basis, and define $V_f=\text{span}\{\gamma_1,\cdots,\gamma_n\}$ to be the
 fine space.
 We obtained the fine solution denoted by $p_h$ by solving
 \begin{eqnarray}
a(p_h,v_h)=F(v_h)-\left\langle g_N,v_h\right\rangle_{\Gamma_{N}} \text{for all }v_h\in V_f \label{fine}
 \end{eqnarray}
  The construction of multiscale basis functions follows two general steps. First, we construct snapshot basis functions in order to build a set of possible modes of the solutions. In the second step, we construct multiscale basis functions with a suitable spectral problem defined in the snapshot space. We take the first few dominated eigenfunctions as basis functions. Using the multiscale basis functions, we obtain a reduced model.
More specifically, once the coarse and fine grids are given, one may construct the multiscale basis functions to approximate the solution
of (\ref{varational}). To obtain the multiscale basis functions, we first define the snapshot space. For each coarse neighborhood $ D_{i}$, define $J_h( D_{i})$ as the set of the fine nodes of $T^{h}$ lying on $\partial D_{i}$ and denote
the its cardinality by $L_i \in \mathbb{N}^{+}$. For each fine-grid node $x_j \in J_h( D_{i})$, we define a fine-grid function $\delta_{j}^{h}$ on $J_h( D_{i})$ as $\delta_{j}^{h}(x_k)=\delta_{j,k}$. Here
 $\delta_{j,k}=1$ if $j=k$ and $\delta_{j,k}=0$ if $j\neq k$. For each $j=1,\cdots, L_i$, we define the snapshot basis functions $\psi_{j}^{(i)}$ ($j=1,\cdots,L_i$) as the solution of the following system\\
\begin{eqnarray}
\begin{aligned}
   -\nabla\cdot\left(\kappa \nabla \psi_{j}^{(i)}\right) &=0 \quad \text { in } D_{i} \\
   \psi_{j}^{(i)} &=\delta_{j}^{h} \quad \text { on } \partial D_{i}.\label{snap_basis}
\end{aligned}
\end{eqnarray}
 The local snapshot space $V_{\text { snap }}^{(i)}$ corresponding to the coarse neighborhood $ D_{i}$  is defined as follows
  $V_{snap}^{(i)}:=$ \text{span}$\{\psi_{j}^{(i)}:j=1,\cdots,L_{i}\}$ and the snapshot space reads $V_{\text {snap}} :=\bigoplus_{i=1}^{N_{\text {in}}} V_{\text {snap}}^{(i)}$, where $N_{\text {in}}$ is the total number of coarse neighborhood.

  In the second step, a dimension reduction is performed on $V_{\text {snap}}$.
  For each $i=1,\cdots, N_{\text {in}}$, we solve the following spectral problem:
  \begin{eqnarray}
      \int_{D_{i}} \kappa \nabla \phi_{j}^{(i)} \cdot \nabla v=\lambda_{j}^{(i)} \int_{D_{i}} \hat{\kappa} \phi_{j}^{(i)} v \quad \forall v \in V_{\text {snap}}^{(i)}, \quad j=1, \ldots, L_{i}\label{eigen}
  \end{eqnarray}
  where $\hat{\kappa} :=\kappa \sum_{i=1}^{N_{i n}} H^{2}\left|\nabla \chi_{i}\right|^{2}$ and $\{\chi_{i}\}_{i=1}^{N_{i n}}$ is a set of partition of unity that solves the following system:
  \begin{eqnarray*}
  \begin{array}
  {rlrl}{-\nabla \cdot\left(\kappa \nabla \chi_{i}\right)} & {=0} & {} & {\text { in } K \subset D_{i}} \\ {\chi_{i}} & {=p_{i}} & {} & {\text { on each } \partial K \text { with } K \subset D_{i}} \\ {\chi_{i}} & {=0} & {} & {\text { on } \partial D_{i}}
  \end{array}
  \end{eqnarray*}
  where $p_i$ is some polynomial functions and we can choose linear functions for simplicity.
  Assume that the eigenvalues obtained from (\ref{eigen}) are arranged in ascending order and we may use the first $1<l_i \leq L_{i}$ (with $l_{i} \in
  \mathbb{N}^{+}$) eigenfunctions (related to the smallest $l_i$ eigenvalues) to
  form the local multiscale space $V_{\text{off}}^{(i)}:=$ snap$\{\chi_{i}\phi_{j}^{(i)}:j=1,\cdots,L_{i}\}$. The mulitiscale space $V_{\text{off}}^{(i)}$ is the direct sum of the local mulitiscale spaces, namely
  $V_{\text {off}} :=\bigoplus_{i=1}^{N_{\text {in}}} V_{\text {off}}^{(i)}$.
  Once the multiscale space $V_{\text {off}}$ is constructed, we can find the
  GMsFEM solution $p_H$ by solving the following equation
 \begin{eqnarray}
a(p_H,v_H)=F(v_H)-(g_N,v_H)_{\Gamma_{N}} \text{for all }v_H\in V_{\text {off}} \label{coarse}
 \end{eqnarray}
 In the numerical examples, we use $L_z$ to denote $L_i$ for $1\leq i\leq N_{\text {in}}$ since we use same $L_i$ for each $i$.
\subsection{Online enrichment}
We will present the constructions of online basis functions \cite{chung2015residual} in this section. \\
After obtaining the multiscale space $V_{\text {off}}$, one may add some online basis functions based on local residuals.\\
Let $p_H \in V_{\text {off}}$ be the solution obtained in (\ref{coarse}). Given a coarse neighborhood $D_i$, we define
$V_i:=H_0^1(D_i)\cap V_{\text {snap}}$ equipped with the norm
$\|v\|_{V_i}^{2}:=\int_{D_i}\kappa |\nabla {v}|^2$. We also define the local residual operator $R_i: V_i\rightarrow \mathbb{R}$ by
\begin{eqnarray}
    \mathcal{R}_{i}\left(v ; p_H\right) :=a(p_H,v)-F(v)+(g_N,v)_{\Gamma_{N}} \label{loc res}
\end{eqnarray}
 The norm of operator $R_i$, denoted by $\|R_i\|_{V_{i}^{*}}$, gives a measure of the quantity of residual.

  Suppose one needs to add one new online basis $\phi$ into the space $V_i$. The analysis in \cite{chung2015residual} suggests that the required online basis $\phi\in V_i$ is the solution to the following equation
  \begin{eqnarray}
      \mathcal{A}(\phi, v)=\mathcal{R}_{i}\left(v ; p_H^{\tau}\right) \quad \forall v \in V_{i}.\label{online}
  \end{eqnarray}
  We refer to $\tau \in \mathbb{N}$  as the level of the enrichment
  and denote the solution of (\ref{coarse}) by $p_H^{\tau}$.
  Remark that $V_{\text {off}}^{0}:=V_{\text{off}}$. Let $\mathcal{I} \subset\left\{1,2, \ldots, N_{i n}\right\}$ be the index set over some non-lapping coarse neighborhoods. For each $i\in \mathcal{I}$, we obtain a online basis $\phi_i\in V_i$ by solving (\ref{online}) and define
$V_{\text {off}}^{\tau+1}=V_{\text {off}}^{\tau} \oplus \operatorname{span}\left\{\phi_{i} : i \in \mathcal{I}\right\}$.
After that, solve (\ref{coarse}) in $V_{\text {off}}^{\tau+1}$.

\section{Postprocessing GMsFEM solution}\label{postprocess}
In order to obtain the GMsFEM with local conservation property, we apply postprocessing technique after obtain $p_H$. The technique was introduced in \cite{bush2014application}. In this section, a review is presented.\\
This approach is composed of two main steps. The first step is solving an auxiliary boundary value problem element by element.  The next step is called downscaling procedure, which is solving similar boundary value problem in each control volumn using the auxiliary solutions obtained in first step. Derivation of local conservation of mass is presented in \ref{conservation}.
\subsection{Constructing a locally conservative flux}
In particular, we obtain a auxiliary solution denoted by $\tilde{p}_{\tau}$,with
\begin{eqnarray}
    \left\{\begin{array}{ll}
-\nabla \cdot\left(\lambda \kappa(x) \nabla \tilde{p}_{\tau}\right)=q & \text { in } \tau \\
-\lambda \kappa(x) \nabla \tilde{p}_{\tau} \cdot n=\tilde{g}_{\tau} & \text { on } \partial \tau
\end{array}\right.\label{auxiliary}
\end{eqnarray}
 \begin{figure}[!htbp]
    \centering
    \subfigure
       { \includegraphics[width=0.45\textwidth]{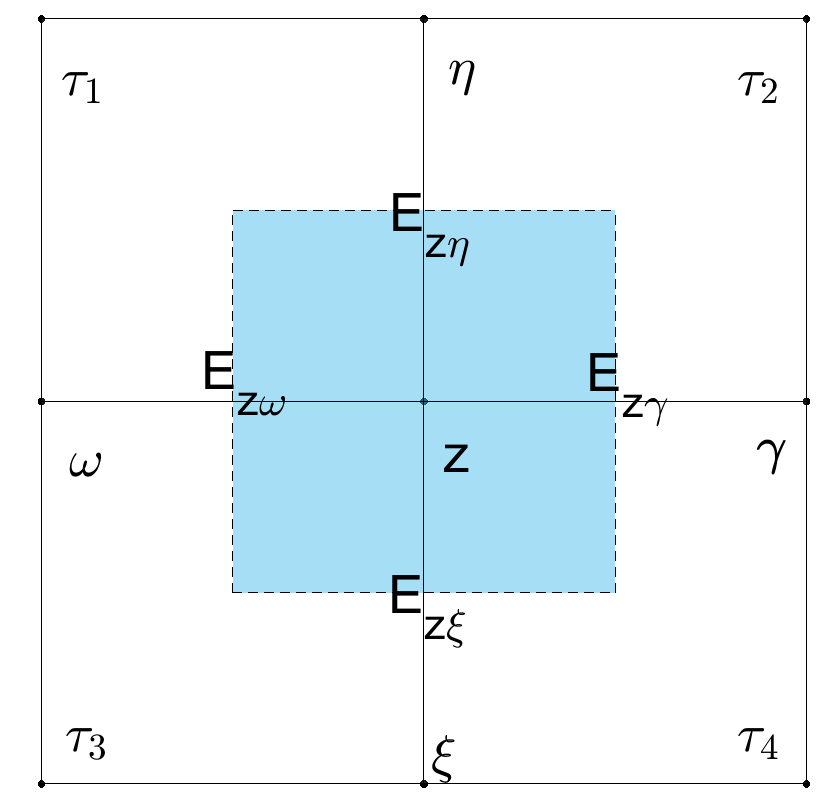}}
        \subfigure
       { \includegraphics[width=0.3\textwidth]{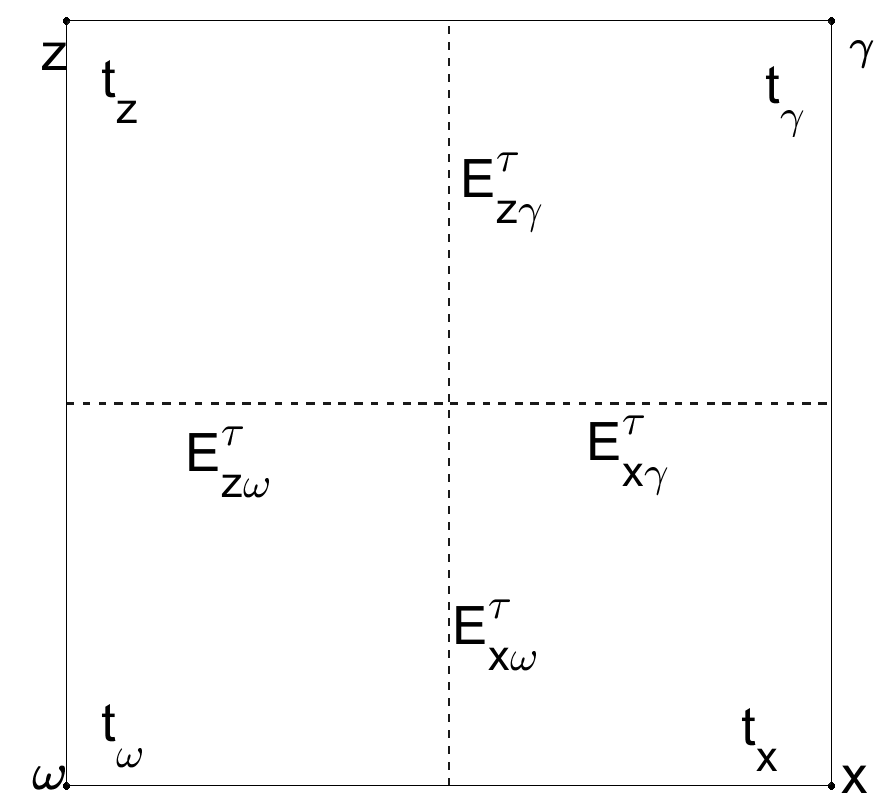}}
    \caption{Left: $C_z$ is the control volumn associated with the vertex z,where $\partial C_z=E_{z\eta}\cup E_{Z\omega}\cup E_{z\xi}\cup E_{z\gamma}$.
    Right: a finite element $\tau$ is divided into four quadrilaterals $t_z,t_{\omega},t_x,t_{\gamma}$.}\label{control volumn}
\end{figure}
Here, we designate $\partial{\tau}=\cup_{\xi\in v(\tau)} E_{\xi}^{\tau}$, where $E_{\xi}^{\tau}=\partial \tau \cap \partial t_{\xi}$ (i.e. half of each element edge containing the vertex $\xi$.) and $v(\tau)$ is the collection of four vertexes of $\tau$.
Furthermore, we set $\tilde{g}_{\tau}$ as piecewise function on $\partial\tau$ such that
\begin{eqnarray*}
    \int_{E_{\xi}^{\tau}} \tilde{g}_{\tau} \mathrm{d} l=F_{\xi, 1}-Q_{\xi, 1}, \quad \text { for } \xi \in v(\tau)
\end{eqnarray*}
where
\begin{eqnarray}
Q_{\xi, 1}=\int_{\tau} \lambda \kappa \nabla p_{H} \cdot \nabla \Phi_{\xi, 1} \mathrm{d} x \quad \text { and } \quad F_{\xi, 1}=\int_{\tau} q \Phi_{\xi, 1} \mathrm{d} x \label{Q and F}
\end{eqnarray}
The existence and uniqueness of the above problem is stated in \cite{bush2014application}.
\ref{auxiliary} implies
\begin{equation*}
    -\int_{\partial \tau} \lambda \kappa \nabla \tilde{p}_{\tau} \cdot \boldsymbol{n} \mathrm{d} l=-\int_{\partial \tau} \lambda \kappa \nabla p \cdot \boldsymbol{n} \mathrm{d} l,
\end{equation*}
which shows that the solution of \ref{auxiliary} recovers the flux of p (i.e. the true pressure solution) averaged over $\partial \tau$, a local conservation property in each element. We use \ref{auxiliary} as a governing principle to derive the processing technique for calculating a locally conservative flux in each control volumn from $p_H$.
 The elemental calculation is based on discretization of $\tau$ into quadrilaterals $t_{\xi}$,i.e., $\tau=\cup_{\xi\in v(\tau)}t_{\xi}$, each of which yields $t_{\xi}=C_{\xi}\cap\tau$,see the right plot of Figure \ref{control volumn}. We set the local solution space as $\mathcal{V}(\tau)=\text{span}\{\Phi_{\xi,1}\}_{\xi\in v(\tau)}$, where $\Phi_{\xi,1}$ is the multiscale basis function corresponding to the vertex $\xi$ . The numerical solution associated with
 \ref{auxiliary} is to find $\tilde{p}_{\tau,h}\in \mathcal{V}(\tau)$ satisfying
 \begin{equation}
     -\int_{\partial t_{\zeta}} \lambda \kappa \nabla \tilde{p}_{\tau, h} \cdot \boldsymbol{n} \mathrm{d} l=\int_{t_{\zeta}} q \mathrm{d} x, \quad \text { for all } \zeta \in v(\tau)\label{part}
 \end{equation}
 The following four equations result from \ref{part}:
 \begin{eqnarray}
     \begin{aligned}
&q_{z \omega}^{\tau}+q_{z\gamma}^{\tau}=Q_{z, 1}-F_{z, 1}+\int_{t_{z}} q \mathrm{d} x\\
&q_{z\gamma}^{\tau}+q_{x \gamma}^{\tau}=Q_{\gamma, 1}-F_{\gamma, 1}+\int_{t_{\gamma}} q \mathrm{d} x\\
&q_{x \omega}^{\tau}+q_{z \omega}^{\tau}=Q_{\omega, 1}-F_{\omega, 1}+\int_{t_{\omega}} q \mathrm{d} x\\
&q_{x \gamma}^{\tau}+q_{x \omega}^{\tau}=Q_{x, 1}-F_{x, 1}+\int_{t_{x}} q \mathrm{d} x \label{system}
\end{aligned}
 \end{eqnarray}
 where
 \begin{eqnarray}
     \begin{aligned}
q_{x\omega}^{\tau} &=-\int_{E_{x\omega}^{\tau}} \lambda \kappa \nabla \tilde{p}_{\tau, h} \cdot \boldsymbol{n} \mathrm{d} l, & & q_{z \gamma}^{\tau}=-\int_{E_{z \gamma}^{\tau}} \lambda \kappa \nabla \tilde{p}_{\tau,h} \cdot \boldsymbol{n} \mathrm{d} l \\
q_{z\omega}^{\tau} &=-\int_{E_{z\omega}^{\tau}} \lambda \kappa \nabla \tilde{p}_{\tau, h} \cdot \boldsymbol{n} \mathrm{d} l, & & q_{x \gamma}^{\tau}=-\int_{E_{x \gamma}^{\tau}} \lambda \kappa \nabla \tilde{p}_{\tau, h} \cdot \boldsymbol{n} \mathrm{d} l
\end{aligned}
 \end{eqnarray}
 and $E_{\xi\eta}^{\tau}=\partial t_{\xi}\cap \partial t_{\eta}$, for $\xi,\eta=\omega,x,\gamma,z$ and $\xi\neq\eta$.
 Since we actually use linear combination of basis to solve solution, in particular,
 $\tilde{p}_{\tau,h}=\Sigma_{\xi\in v(\tau)} u_{\xi}\Phi_{xi}$ with unknown coefficients $u_{\xi}$,
 \ref{system} can be written in the form of
 $\tilde{A}\tilde{u}=\tilde{f}$ where

 \begin{equation*}
 \tilde{\boldsymbol{A}}_{\zeta \eta}=-\int_{E_{\mathrm{fp}}} \lambda \kappa \nabla \Phi_{\eta,1}+\boldsymbol{n} \mathrm{d} l \quad \text { and } \quad \tilde{\mathrm{f}}_{\zeta}=\int_{t_{\mathrm{\xi}}} q \mathrm{d} x-\int_{E_{\xi}^{\tau}} \tilde{g}_{\tau} \mathrm{d} l
 \end{equation*}
One should note that when $\tau$ is adjacent to $\Gamma_N$, $g_N$ should be taken in account in computing $\tilde{g}_{\tau}$.\\
Since the system actually has smaller dimension than 4, we may add a constant to one entry in $\tilde{A}$ to remove the singularity. The fact that $u$ is not unique is irrelevant since the desired solution is flux as governed by $q_{\xi\eta}^{\tau}$ which is unique.\\
\ref{system} implies that $\tilde{v}_h$ derived from $\tilde{p}_h$ satisfy the desired local conservation property.
\subsection{Downscale procedure}
After the postprocessing in the section 4.1, we have
\begin{equation*}
    \int_{\partial C_{\tau}} \tilde{v}_{h} \cdot \boldsymbol{n} \mathrm{d} l=\int_{\mathcal{C}_{z}} q \mathrm{d} \boldsymbol{x}, \quad \text { for all } C_{z}
\end{equation*}
which can be thought of a statement of compatibility condition in $C_z$.
We can proceed with formulating a boundary problem as follows,
\begin{eqnarray}
    \left\{\begin{array}{l}
-\nabla \cdot\left(\lambda \kappa(\boldsymbol{x}) \nabla \tilde{p}_{G}\right)=q \quad \text { in } C_{z} \\
-\lambda \kappa(\boldsymbol{x}) \nabla \widetilde{p}_{G} \cdot \boldsymbol{n}=\tilde{v}_{h} \cdot \boldsymbol{n} \text { on } \partial C_{z}
\end{array}\right.\label{downscale}
\end{eqnarray}
 Here $\tilde{v}_{h}=\Sigma_{\tau,\tau\cap C_z\neq \emptyset}-\lambda \kappa(x)\nabla\tilde{p}_{\tau,h}$ that is evaluated pointwise on segments of $\partial C_z$ that belongs to $\tau$.\\
  For example, for control volume $C_z$ corresponding to vertex $z$, we obtain $\tilde{v}_{h}$ as follows. One may refer to the left in figure \ref{control volumn}.
  \begin{eqnarray}
      \begin{aligned}
\int_{\partial C_{z}} \tilde{v}_{h} \cdot n d l &=\int_{E_{z\eta}} \tilde{v}_{h} \cdot n d l+\int_{E_{z\omega}} \tilde{v}_{h} \cdot n d l+\int_{E_{z\xi}} \tilde{v}_{h} \cdot n d l+\int_{E_{z\gamma}} \tilde{v}_{h} \cdot n d l \\
&=\left(q_{z\eta}^{\tau_{1}}+q_{z\eta}^{\tau_{2}}\right)+\left(q_{z \omega}^{\tau_{1}}+q_{z\omega}^{\tau_{3}}\right)+\left(q_{z\xi}^{\tau_{3}}+q_{z\xi}^{\tau_{4}}\right)+\left(q_{z \gamma}^{\tau_{4}}+q_{z \gamma}^{\tau_{2}}\right)\\
&=\left(q_{z\eta}^{\tau_{1}}+q_{z\omega}^{\tau_{1}}\right)+\left(q_{z\eta}^{\tau_{2}}+q_{z\gamma}^{\tau_{2}}\right)+\left(q_{z\omega}^{\tau_{3}}+q_{z_{z\xi}}^{\tau_{3}}\right)+\left(q_{z\xi}^{\tau_{4}}+q_{z\gamma}^{\tau_{4}}\right)\\
&=\sum_{j=1}^{4}\left(Q_{z, 1, j}-F_{z, 1, j}\right)+\sum_{j=1}^{4} \int_{t_{z, j}} q \mathrm{d} \mathbf{x}\\
&=\int_{ C_{z}} q \mathrm{d} \mathbf{x}. \label{conservation}
\end{aligned}
 \end{eqnarray}
 where $Q_{z, 1, j},F_{z, 1, j}$ are integrals (refer \ref{Q and F}) in domain $\tau_j$ for corresponding $j$. $\sum_{j=1}^{4}\left(Q_{z, 1, j}-F_{z, 1, j}\right)=0$ is derived by \ref{coarse}.\\
This calculation actually proves the local conservation of $\tilde{v}_{h}$. So the satisfies compatibility condition of \ref{downscale} guarantees the existence of the corresponding solution. Similarly, since our interest only lies in $-\lambda \kappa(x)\nabla\tilde{p}_{C_z}$ in \ref{downscale}, the nonuniqueness of the solution is of no concern.\\

\section{Numerical results}
In this section, we consider four kinds of permeability coefficients which are represented in figure \ref{medium}. $\kappa_1$ and $\kappa_4$ are extracted from the tenth SPE comparative solution project (SPE10), which is commonly used as benchmark permeability field to assess upscaling and multiscale methods. The most distinguishable characteristic of the model is that some layers are highly heterogeneous and contains long channels. Here, $\kappa_1$ is the last layer of the SPE10 dataset while $\kappa_4$ comes from the 36-th layer. It is evident that $\kappa_1$ represents high heterogeneity and both two contains some visible channels. In terms of $\kappa_2$, it is deterministic, high-contrast coefficient with abrupt transitions between regions of low and high permeability. For $\kappa_3$, it comes from fractured porous media, which is characterized by complex fracture distribution and high contrast. Consequently, four examples of permeability exhibit high-contrast features, which can make solving (\ref{model}) a demanding task.

Since the construction of multiscale space is based on the single-phase flow, i.e. choosing $\lambda(S)=1$ in (\ref{model}), it is reasonable to consider the efficiency of our approximation space within the context of single-phase and further estimate the effect on the two-phase model. Both the two models are solved in the domain $\Omega=[0,1]\times[0,1]$.

\subsection{Single-phase flow}
In the single-phase model, we solve the pressure equation (\ref{model}) with $\lambda=1$ and postprocess the velocity field with the technique introduced in section \ref{postprocess}. For boundary condition, we set Dirichlet boundary condition $p_{D}=1$ on the left edge and $p_{D}=0$ on the right edge of the domain. Besides, we set $g_{N}=0$ on $\Gamma_{N}$, i.e. zero Neumann boundary condition for bottom and top edge. We assume there is no external force so we take $q=0$.  The size of the permeability coefficient is $100*100$ and $200*200$ for $\kappa_2$ and $\kappa_3$ accordingly. To estimate our method GMsFEM, we compare four method, the standard finite element method and GMsFEM with different combination of multiscale basis functions. For GMsFEM, we set the coarse mesh size to be $10*10$. In other words, there are $20*20$ coarse elements in the whole domain. As is shown in \ref{tab:v_error}, there is significant error decay in both cases, where online enrichment contributes much more compared to the offline enrichment. As for the notation, we use $L_z=a+b$ to denote the case where a offline basis followed by b online basis are used in each local neighborhood. In particular, with $\kappa_2$, we can see a sharp decrease from the initial case with big error to a relatively low error when we add the both offline and online basis to the case $2+1$, which is even slightly lower than case $5+0$. In other words, the information contained in online basis functions results in bigger help than that in offline basis. This is due to the global construction of online basis functions while the offline space is constructed locally. As to the other case with $\kappa_3$, we can see the decay is less pronounced, however, there is still evident improvement in the accuracy. It is similar here that we can use online basis functions to obtain satisfying results with smaller dimension of multiscale space.

\begin{figure}[!htbp]
    \includegraphics[width=5in]{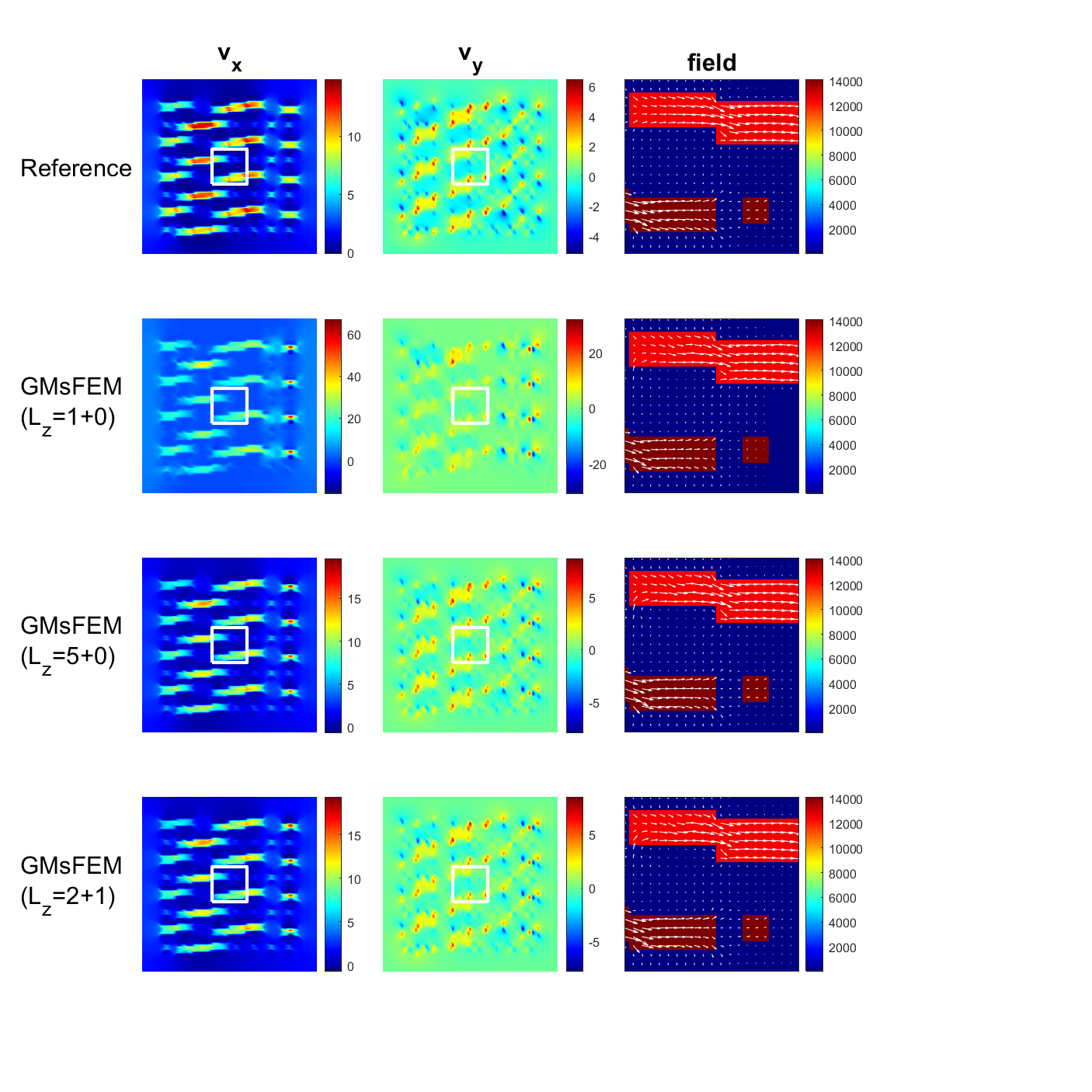}\caption{$\kappa_2$. Velocity computed using four methods. The first two columns (horizontal the first and vertical the second ) exhibits velocity profile on the whole domain with the reference on the first row and last three corresponding to $L_z=1$, $L_z=5$ and $L_z=2+1$ respectively. The last column shows the velocity  with permeability in the selected region.}
    \label{v2}
\end{figure}

\begin{figure}[!htbp]
    \includegraphics[width=5in]{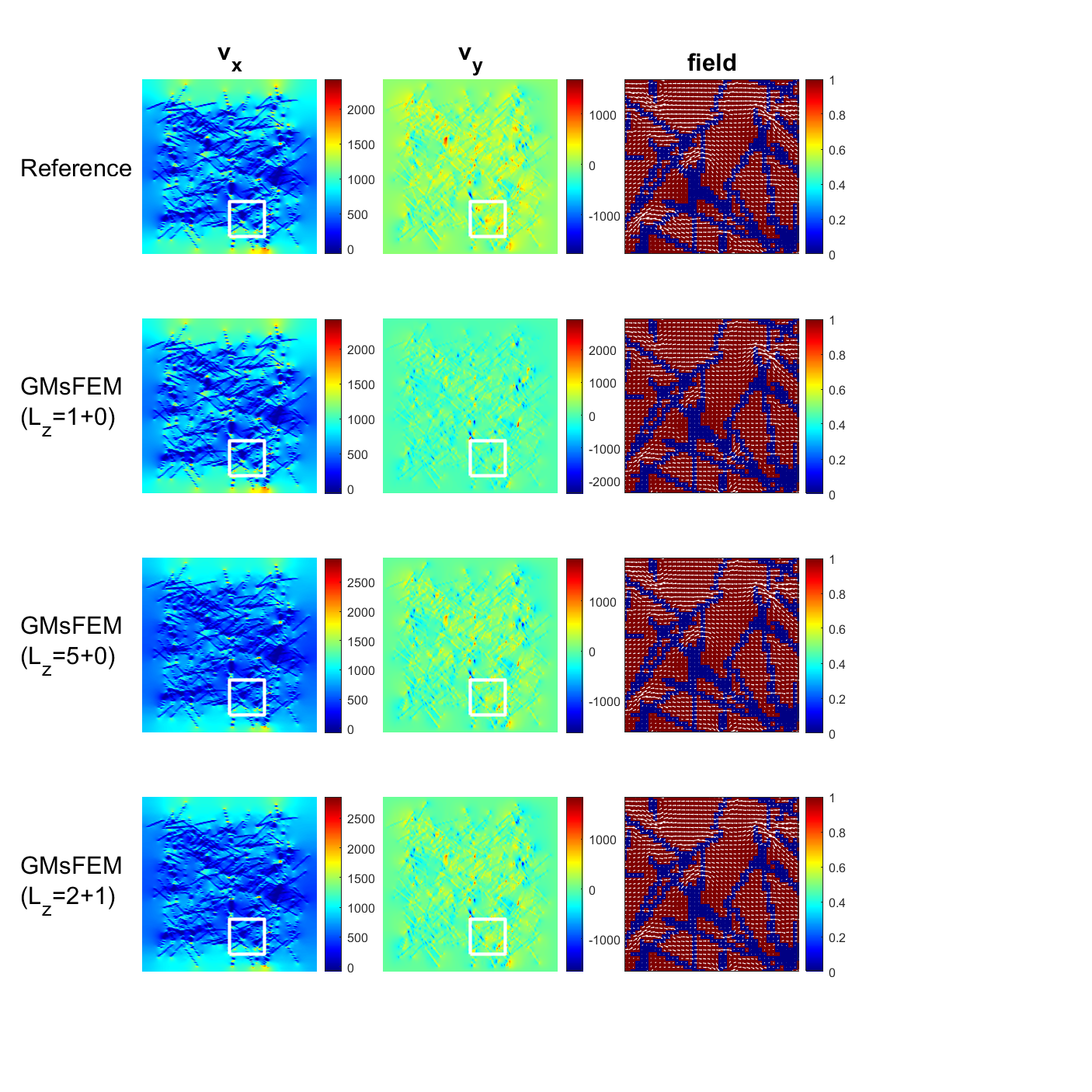}\caption{$\kappa_3$. Velocity computed using four methods. The first two columns (horizontal the first and vertical the second ) exhibits velocity profile on the whole domain with the reference on the first row and last three corresponding to $L_z=1$, $L_z=5$ and $L_z=2+1$ respectively. The last column shows the velocity  with permeability in the selected region.}
    \label{v3}
\end{figure}

To better present the approximation of the velocity field, which is actually used in the further two-phase flow, we plot figure \ref{v2} and \ref{v3} for $\kappa_2$ and $\kappa_3$. In both cases, we show the horizontal and vertical components of velocity. Since the velocity is highly related with the permeability, we exhibit the velocity field under some region with big contrast in permeabilty coefficients, which are also shown as background. We can observe significant dismatch between the initial case and reference while in the case $L_z=2+1$ and $5+0$, the accuracy improvement is apparent especially in the selected region, where the permeability changes rapidly. Specifically, in \ref{v2}, there is a few flows with opposite direction for $L_z=1+0$ compared with the reference, while the difference is less noticeable in the latter two cases.

\begin{table}[]
    \centering
    \begin{tabular}{c|c|c}
    \hline\hline
                  & $\kappa_2$ & $\kappa_3$ \\
                  \hline
       $L_z$=1+0  &  1.30  & 0.21\\
       $L_z$=5+0  &  0.16  & 0.09\\
       $L_z$=2+1  &  0.14  & 0.08\\
    \end{tabular}
    \caption{Relative error for the velocity field. The $L^2$ error of the velocity is computed for different choices of GMsFEM compared with reference solution obtained in standard finite element method.}
    \label{tab:v_error}
\end{table}
\subsection{Two-phase flow}
In solving \ref{transport}, we use the quadratic relative permeability curves $\kappa_{rw}=S^2$ and $\kappa_{ro}=(1-S)^2$, along with $\mu_w=1$ and $\mu_o=5$ for the water viscosities. The domain $\Omega=[0,1]\times[0,1]$. For the initial condition, the value at the left edge is set as $S = 1$ and we assume $S(x, 0) = 0$ elsewhere.  In practical, we construct the multiscale basis functions within the context of single-phase flow model and apply the resulted approximation space to the interested two-phase problem without updating basis function. In other words, we can precompute the bases as preparation before the simulation, which is efficient compared to the case when we need to repeat the computation for different cases.

\begin{figure}[!htbp]
    \includegraphics[width=5in]{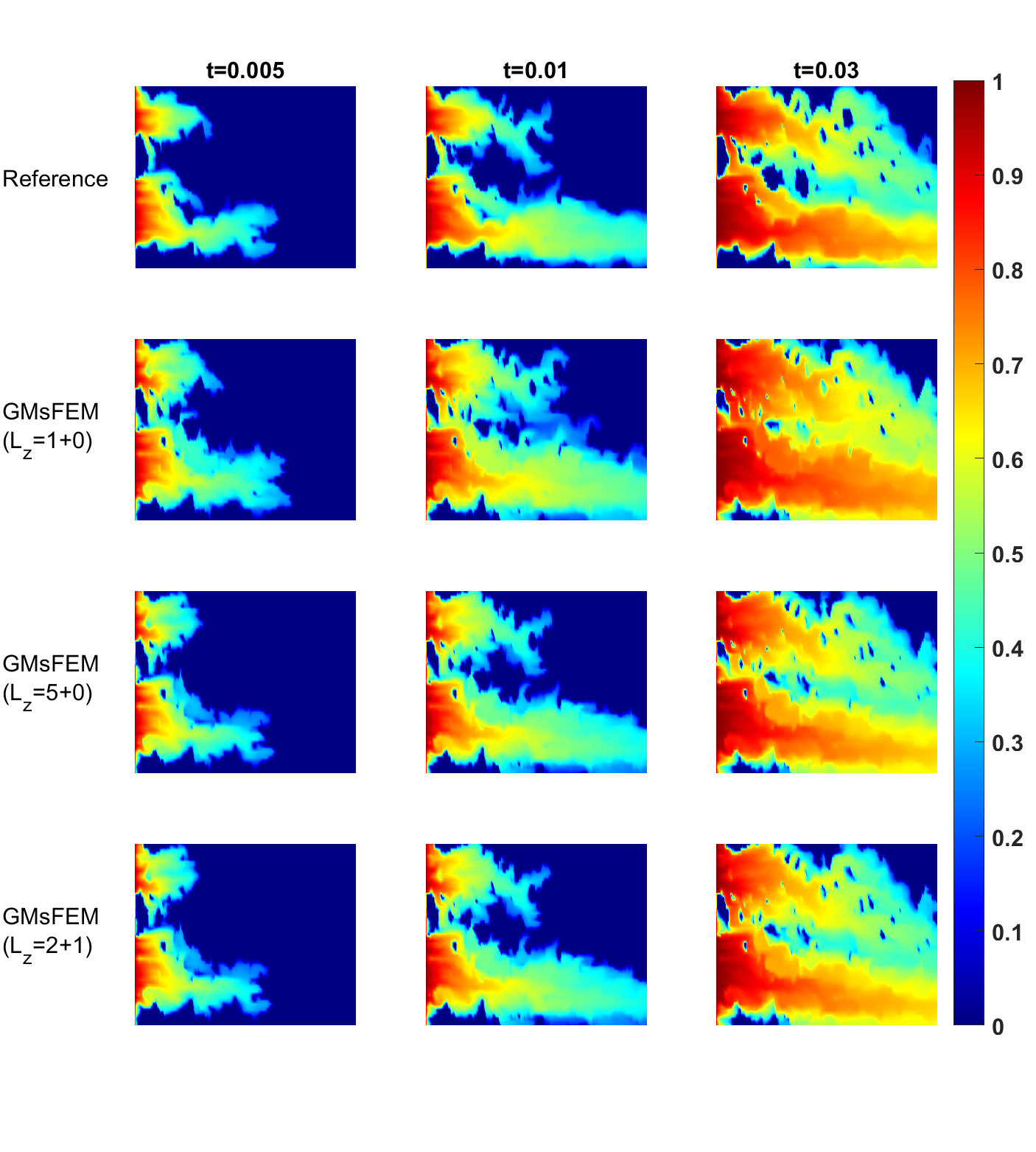}

    \caption{$\kappa_1$. The reference saturation is shown on the top row at three different time levels. The second and third rows are $L_2=1$, $L_z=5$ respectively.The last row is using 2 offline basis and enriched by one online basis for each coarse neighborhood, which is denoted by $L_z=2+1$.}
    \label{k1}
\end{figure}

\begin{figure}[!htbp]
    \includegraphics[width=5in]{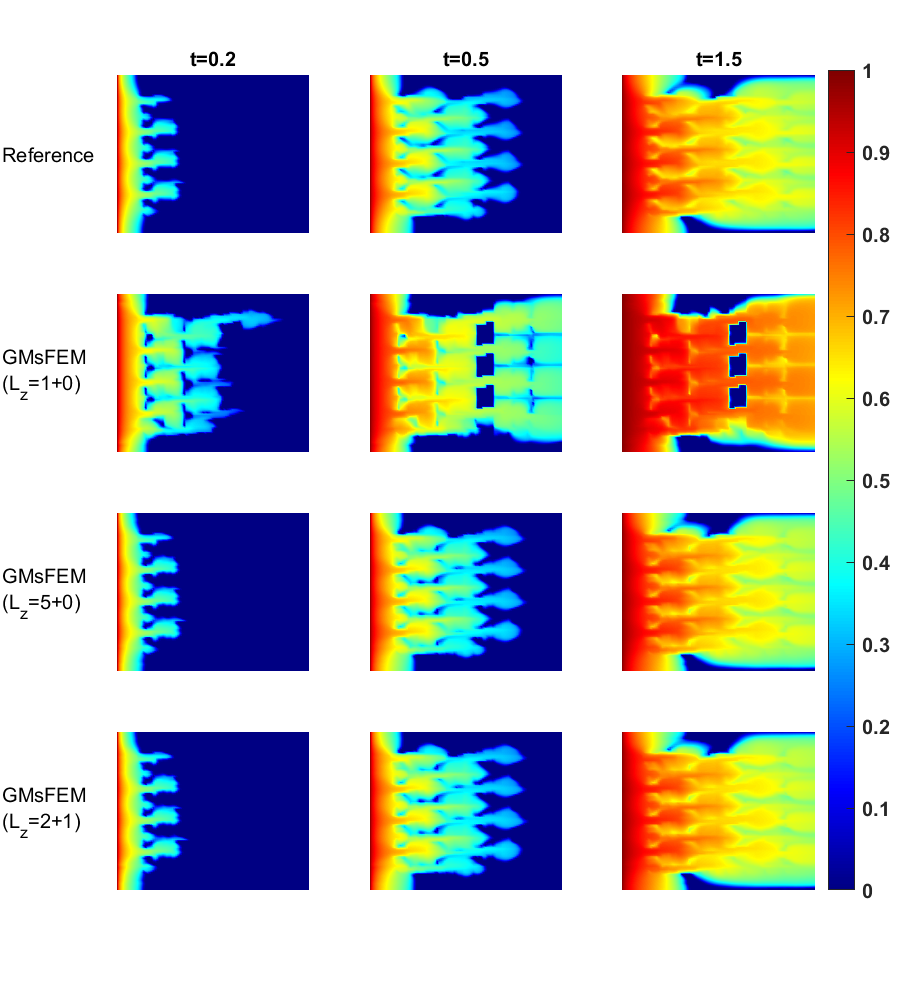}\caption{$\kappa_2$. The reference saturation is shown on the top row at three different time levels. The second through fourth rows are $L_z=1$, $L_z=5$ and $L_z=2+1$ respectively.}
    \label{k2}
\end{figure}

\begin{figure}[!htbp]
    \includegraphics[width=5in]{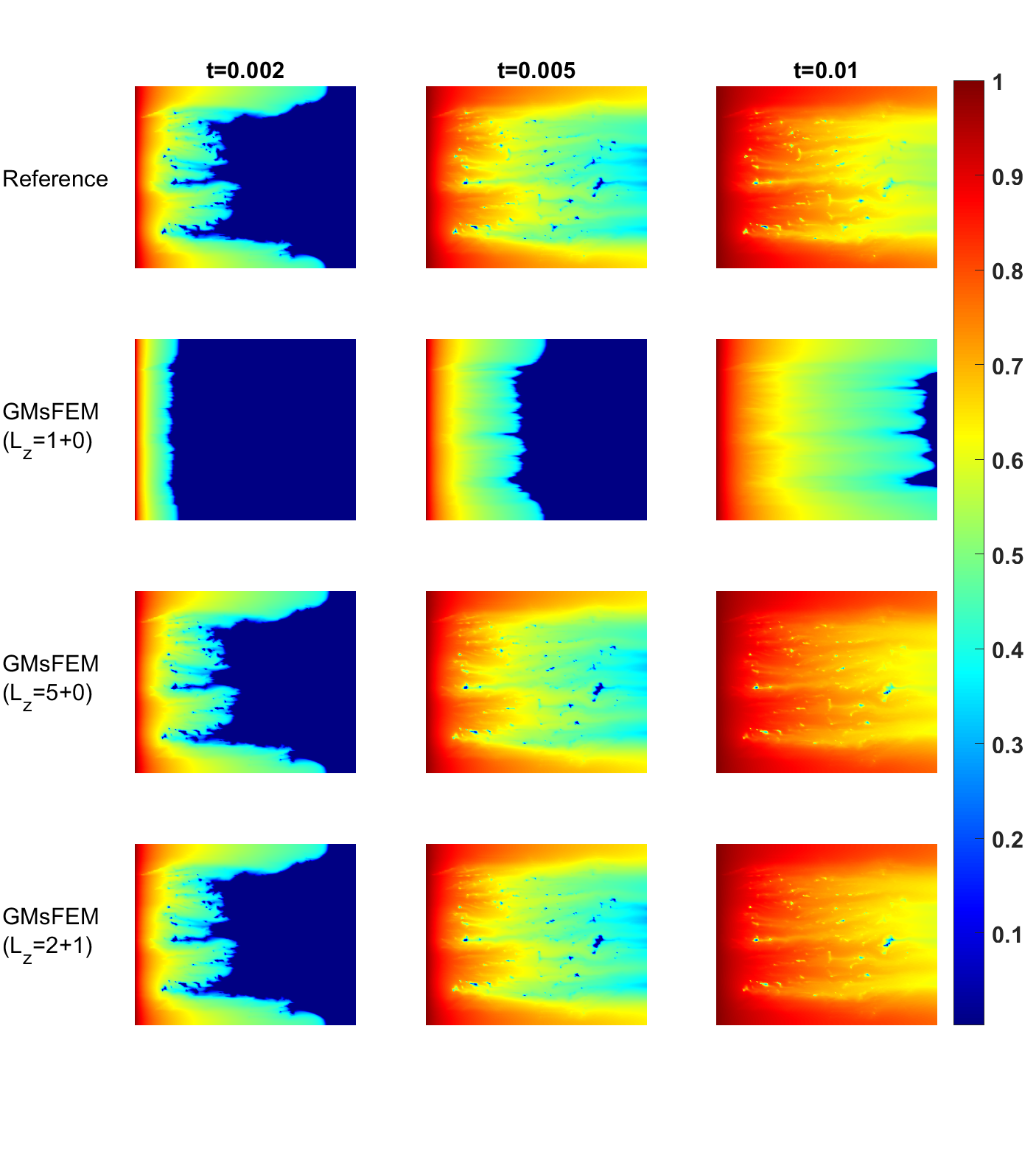}\caption{$\kappa_3$. The reference saturation is shown on the top row at three different time levels. The second through fourth rows are $L_z=1$, $L_z=5$ and $L_z=2+1$ respectively.}
    \label{k3}
\end{figure}

\begin{figure}[!htbp]
    \includegraphics[width=5in]{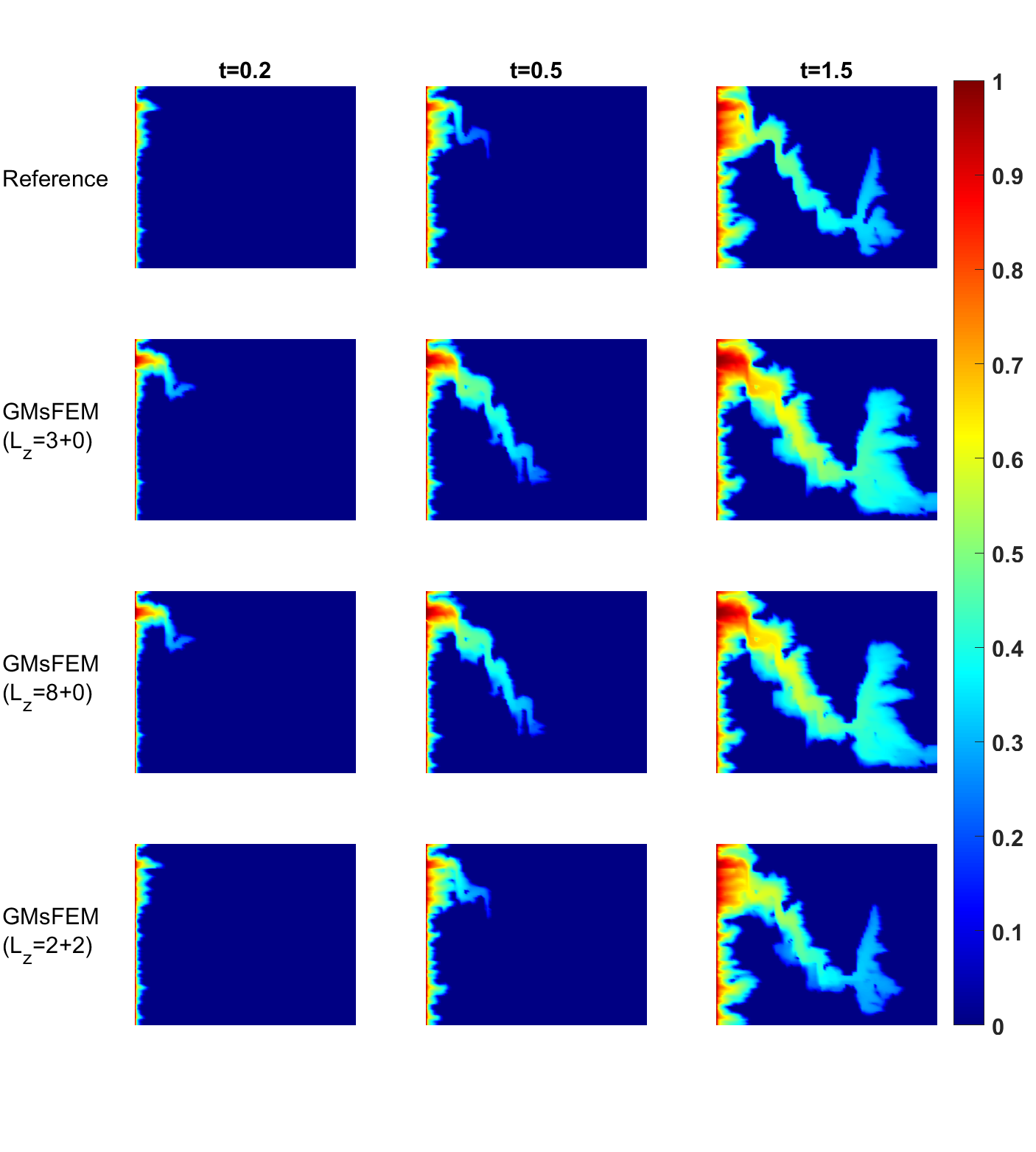}\caption{$\kappa_4$. The reference saturation is shown on the top row at three different time levels. The second through fourth rows are $L_z=3$, $L_z=8$ and $L_z=2+1$ respectively.}
    \label{k4}
\end{figure}

\begin{figure}[!htbp]
\subfigure[$\kappa_1$]
{\includegraphics[width=0.48\textwidth]{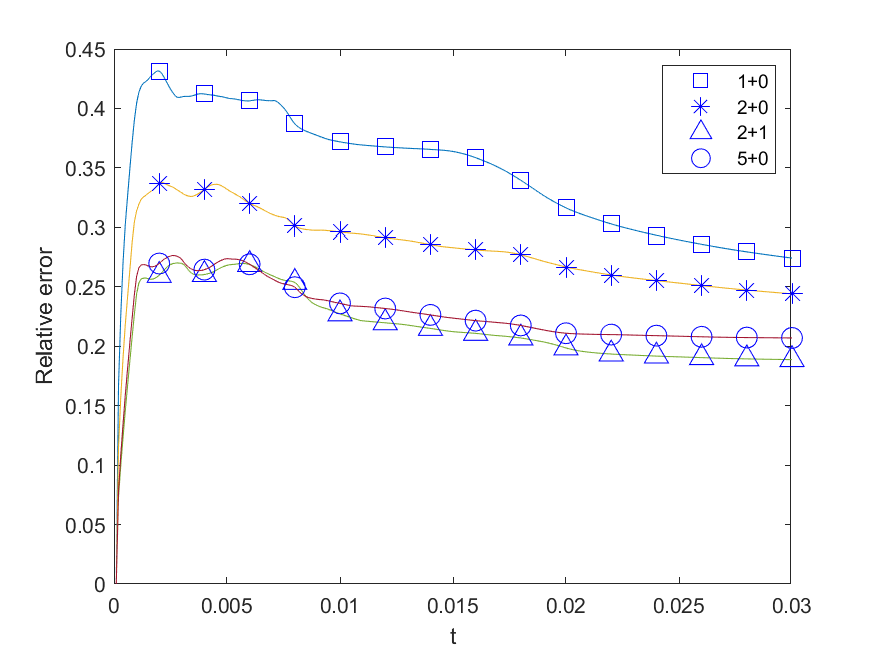}}
\subfigure[$\kappa_2$]
{\includegraphics[width=0.48\textwidth]{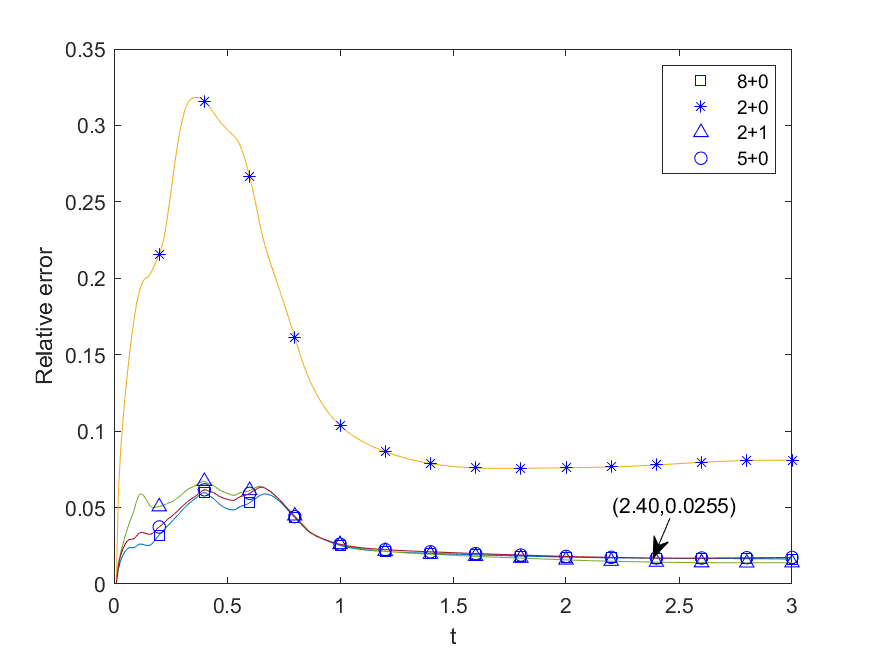}}
\subfigure[$\kappa_3$]
{\includegraphics[width=0.48\textwidth]{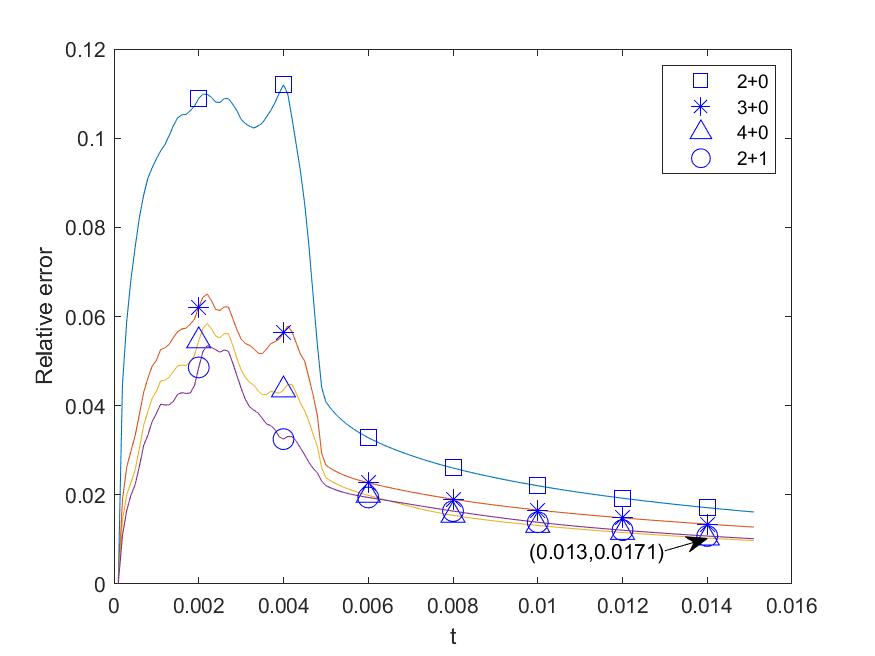}}
\subfigure[$\kappa_4$]
{\includegraphics[width=0.48\textwidth]{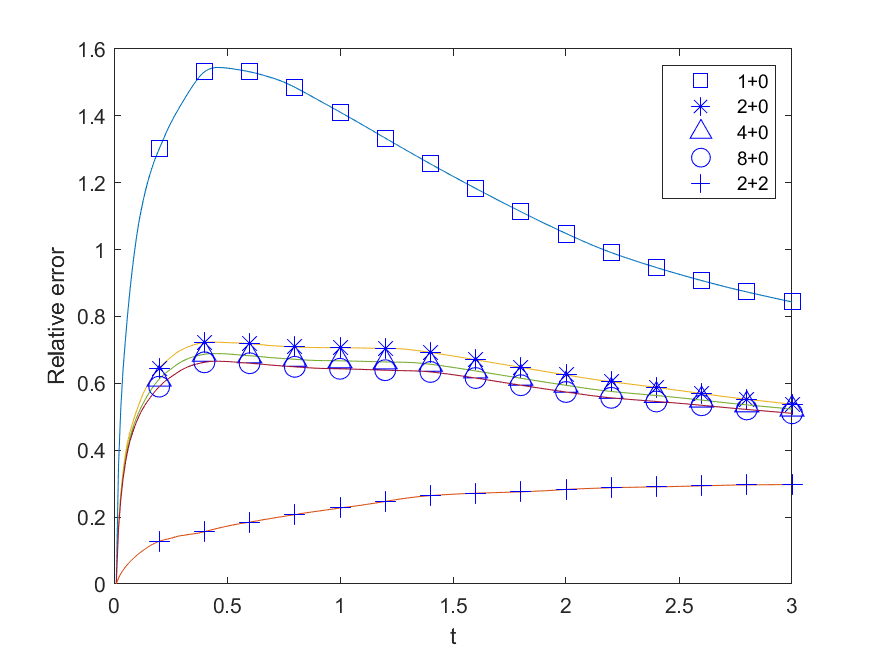}}
\caption{Comparison of the $L^2$ error of the saturation for $\kappa_1$ (upleft) ,$\kappa_2$ (upright) ,$\kappa_3$ (downleft) ,$\kappa_4$ (downright)as a function of time .}\label{error}
\end{figure}

For better visual comparison, we present the saturation for three different time levels in figure \ref{k1},\ref{k2},\ref{k3} and \ref{k4}. From figure \ref{k2} and figure \ref{k3}, significant difference from reference saturation can observed when $L_z=1$ while $L_z=2+1$ are relatively indistinguishable from reference. For case $\kappa_1$ and $\kappa_4$, the improvements are less significant yet pronounced since there are few noticeable differences between last row and reference row.
In other words, online basis functions efficiently improves accuracy compared with offline case. In figure \ref{k4}, $L_z=2+1$ is a better approximation of reference even than $L_z=8+0$. As we can verify it from figure \ref{error}, there are sharp decreases by enriching multiscale space from intial state especially with $\kappa_2$ and $\kappa_3$ compared to the other two cases, which is consistent with the previous dynamics of saturation.
Furthermore, the relative errors are improved by increasing the number of $L_z$ up to a certain threshold, and then the reduction is minimal as more functions are added. In particular, in figure (d), when we double the number of offline basis functions from the very beginning, i.e. single in each local neighborhood, the error reduction is evident however improvement is indistinguishable from $2+0$ to $8+0$, which means very limited reduction can be expected by further increasing offline basis functions. At the same time, adding a few online basis functions will notably increase the accuracy since the error in $L_z=2+2$ is even lower than the case $L_z=8$, which shows the power of incorporating global information inherited in the online basis functions. For $\kappa_1$ and $\kappa_3$, the reduction resulted by enrichment is relatively steady while in $\kappa_2$, it is easier to reach a threshold. This is due to higher heterogeneity in $\kappa_1$ compared with $\kappa_2$. However, all the above four cases combined with the single-phase case share the same conclusion that online enrichment offer us better accuracy with relatively lower cost. Therefore, it is efficient to compute residual-driven online basis functions for the sake of increasing accuracy. \\
\section{Conclusion}
In this paper, we consider a conservation GMsFEM for treating the coupled pressure-convection-diffusion system in the context of the two-phase flow model. An advantage of the proposed method is that local conservative and accurate velocity field can be obtained, which means we combined two main procedures, postprocessing and online enrichment. The effect can be verified in the numerical results. In the future, we can work on more computational-efficient ways to achieve the goal.

\bibliographystyle{IEEEtran}
\bibliography{reference}
\end{document}